\documentclass[12pt]{paper}

\usepackage{amscd}

\usepackage[headings]{fullpage}
\usepackage{mathrsfs}
\usepackage{amsmath,amsfonts,amssymb,amsthm}
\usepackage[all]{xy}
\usepackage{amscd}

\usepackage[dvips]{graphics,color}
\usepackage{hyperref}

\def\noi{\noindent}

\def\pf{\noi{\bf Proof.\ \,}}

\def\eop{{$\square$}}

\def\NN{{\mathbb N}}

\def\ZZ{{\mathbb Z}}
\def\ma{{\mathfrak a}}

\def\mm{{\mathfrak m}}

\def\mp{{\mathfrak p}}
\def\mq{{\mathfrak q}}

\def\Supp{{\mathrm{Supp}}}
\def\Ass{{\mathrm{Ass}}}

\def\dim{{\mathrm{dim}}}
\def\Spec{{\mathrm{Spec}}}
\def\Supp{{\mathrm{Supp}}}

\def\Hom{{\mathrm{Hom}}}
\def\Hhom{{\mathrm{\mathbf{Hom}}}}

\def\Cone{{\mathrm{Cone}}}

\def\perf{{\mathrm{perf}}}
\def\fg{{\mathrm{fg}}}

\def\Ch{{\mathrm{Ch}}}
\def\K{{\mathrm{K}}}

\def\D{{\mathrm{D}}}

\def\ann{{\mathrm{ann}}}

\def\al{ {\cal A}  }

\def\dl{ {\cal D}  }

\def\pl{ {\cal P}  }
\def\tl{ {\cal T}  }

\def\a{\alpha}

\def\t{\tau}

\begin{document}

\newtheorem{thm}{Theorem}[section]
\newtheorem{prop}[thm]{Proposition}
\newtheorem{lem}[thm]{Lemma}
\newtheorem{rem}[thm]{Remark}
\newtheorem{coro}[thm]{Corollary}
\newtheorem{conj}[thm]{Conjecture}
\newtheorem{de}[thm]{Definition}
\newtheorem{hyp}[thm]{Hypothesis}

\newtheorem{nota}[thm]{Notation}
\newtheorem{ex}[thm]{Example}
\newtheorem{proc}[thm]{Procedure}

\newtheorem{exer}[thm]{Exercise}

\newtheorem{prob}[thm]{Problem}

\begin{center}

{\bf \ A classification of nullity classes in the derived category of a ring}
\medskip

Yong Liu \footnote{Email address is \emph{yongliue@gmail.com}} and
Don Stanley\footnote{Email address is
\emph{Donald.Stanley@uregina.ca}}

Department of Mathematics and Statistics

University of Regina

\end{center}

\centerline{version  \today  }

\begin{abstract}
For a commutative Noetherian ring $R$ with finite Krull dimension, we study the nullity classes in $\D^c_\fg(R)$, the full triangulated subcategory $\D^c_\fg(R)$ of the derived category $D(R)$ consisting of objects which can be represented by cofibrant objects with each degree finitely generated. In the light of perversity functions over the prime spectrum $\Spec R$, we prove that there is a complete invariant of nullity classes thus that of aisles (or equivalently, $t$-structures) in $\D^c_\fg(R)$.
\end{abstract}

\medskip

\tableofcontents

\newpage

\section{Introduction}



The classification of various types of algebraic or geometric objects, such as finite simple groups and differentiable manifolds of dimension 4, is a very fundamental problem in many different areas of mathematics. Finding invariants is considered efficient to distinguish them, which in fact establishes a bijection from the collection of target objects to that of well-studied objects (such as finitely generated abelian groups or prime spectrum of a ring). Our research were concerned with the classification of subcategories. These problems naturally arose in homotopy theory (see Hopkins~\cite{Hopkins}) and have heavily influenced algebraic geometry (see Thomason~\cite{Thomason} and Neeman~\cite{Neeman}) and modular representation theory (see Benson, Iyengar and Krause~\cite{Benson}) since then. Explicitly, Gabriel~\cite{Gabriel} for example, showed that in the category of finitely generated modules over a commutative noetherian ring $R$, the Serre subcategories are one-to-one correspondent to the specialization closed subsets of the prime spectrum $\Spec R$. Neeman in~\cite{Neeman} and~\cite{Neeman2} showed that the (co)localizing subcategories of the derived category $D(R)$ are bijectively correspondent to the subsets of $\Spec R$. More recently in the context of tensor triangulated category, the radical thick tensor ideals are classified by the closed subsets (thanks to Hochster dual~\cite{Hochster}) of Balmer's spectrum of the tensor triangulated category, see Balmer~\cite{Balmer}.

In this paper, we are particularly interested in the classification of $t$-structures. The notion of $t$-structures was introduced in the work of Beilinson, Bernstein and Deligne~\cite{BBD}, which commenced the theory from a simple idea of truncation of complexes and connected an abelian category and its derived category in a homological way. Later on, Keller and Vossieck~\cite{Keller} pointed out that $t$-structures are bijectively correspondent to the subcategories aisles (see Theorem~\ref{t-structure}). Nullity classes a generalization of aisles, were introduced by Stanley~\cite{Stan}, who gave a classification of nullity classes (particularly, aisles thus $t$-structures) in the derived categories $\D_\perf(R)$ of perfect complexes and also in $\D^b_\fg(R)$ of bounded derived category whose objects have finitely generated homologies. We generalized his result in the current paper. Recently, translated into our language, Tarrio, Lopez and Saorin~\cite{Alonso2} gave a classification of compactly generated aisles in $\D^+_\fg(R)$ for a commutative Noetherian ring $R$, also via the perversity functions over $\Spec R$.


In our classification of nullity classes, the cellular tower obstruction of objects in $\D(R)$ and their homology supports were used. We strengthened the assumption on the category $\D_\fg(R)$ to $\D^c_\fg(R)$ in which every object has a cofibrant representative with each degree finitely generated, thanks to the model category theory (see Dwyer and Spalinski~\cite{dwyerspalinski} and Hovey~\cite{Hovey}, for example). The paper is organized as follows. In Section~\ref{suppass} and Section~\ref{atn}, we introduce the basic notions such as nullity classes and support of modules. Then in Section~\ref{kca}, we show as Lemma~\ref{genkilling} that to construct an object $M$ in $\D(R)$ it is sufficient to use the quotient $R/\mp$ by primes $\mp$ in the support of $M$. In Section~\ref{key}, a technical lemma is proved as Lemma~\ref{crucial}, which roughly shows the converse statement of Lemma~\ref{genkilling} provided that the ring $R$ has finite Krull dimension. Finally, the complete invariant of nullity classes in $\D^c_\fg(R)$ is given in the last section, also it gives a complete invariant of aisles or $t$-structures in the same category.


Through out the paper, we assume that the ring $R$ is commutative Noetherian with identity.

\section{Support and associated primes}\label{suppandass}\label{suppass}

Let $\D(R)$ be the derived category of a ring $R$. For the definitions and various properties of $\D(R)$ one can see~\cite{weib} for example. We will adopt the following notations for various full triangulated subcategories of $\D(R)$:\\

\noindent(1) The category $\D^\ast_\fg(R)$ ($\ast=b,+,-,\text{blank}$, resp.) consists of (bounded, bounded below, bounded above, blank, resp.)
chain complexes of $R$-modules whose homologies are finitely
generated;

\noindent(2) The category $\D_\perf(R)$ consists of
bounded complexes of finitely generated projective $R$-modules, i.e. \textit{compact} objects, the Hom functors represented by such objects commute with coproducts;

\noindent(3) The category $\D^c_\fg(R)$ consists of objects which can be represented by cofibrant objects finitely generated in each degree. See Proposition~\ref{weakcpt} for detail.\\

We denote by $\Spec R$ the prime spectrum of $R$, i.e. the set of all prime ideals of $R$. For each $R$-module $M$, we can consider two useful subsets of
$\Spec(R)$. Let $\ann^R(x)=\{r\in R~|~rx=0\}$ be the annihilator of $x\in M$ in $R$ and write $M_\mp=M\otimes R_\mp$ for abbreviation. The set of \textit{associated primes} of $M$ is defined as $\Ass M=\{\mp\in
\Spec(R)~|~\mp=\ann^R(x)$ for some $x\in M\}$, and the \textit{support} of $M$ is defined as $\Supp M=\{\mp\in \Spec(R)~|~
M_\mp\neq 0\}$.

Similarly, we can define the set of associated primes and the support of an object $M\in\D(R)$ by setting $\Ass
M=\bigcup_i\Ass H_i(M)$ and $\Supp M=\{\mp\in\Spec R~|~ M\otimes R_\mp\neq 0\}=\bigcup_i \Supp H_i(M)$. Then the following properties about $R$-modules can be generalized to the case of chain complexes of $R$-modules without difficulty.

Let $U\subseteq \Spec(R)$ be a subset. Define $\overline{U}=\{\mp\in
\Spec R~|~\mp\supseteq \mq$ for some $\mq\in U\}$ to be the
\textit{closure} of $U$ \textit{under specialization}. Particularly, we also denote by $V(\mp)$ for $\overline{U}$ when $U=\{\mp\}$.

\begin{lem}~\label{commalg}
Let $A_i, A, B$ and $C$ be $R$-modules and $0\rightarrow A\rightarrow B\rightarrow C\rightarrow 0$ a short exact sequence. Then

\noindent(1) $\Ass\,A\subset \Ass\,B\subset \Ass\,A\cup \Ass\,C$;

\noindent(2) $\Supp\, B=\Supp\, A\cup \Supp\, C$;

\noindent(3) $\Ass\, C\subset\overline{\Ass\,B}=\Supp\, B$;

\noindent(4) $\Supp\bigoplus_iA_i=\bigcup_i\Supp \,A_i$.

\end{lem}

\pf These results can be found in Chapter 3 of~\cite{Atiyah}, Chapter II and IV of~\cite{Bourbaki}.~~~\eop

\begin{ex}\label{suppkp}
Let $\mp$ be a prime ideal and $k(\mp)$ the corresponding residue field. Take $\mq=\ann^R(\frac{r}{s}+\mp R_\mp)\in \Ass\, k(\mp)$, where $r\in R$ and $s\notin\mp$. Then for every $x\in\mq$ we have $x(\frac{r}{s}+\mp R_\mp)=\mp R_\mp$ thus $\frac{x}{1}\in\mp R_\mp$ since $k(\mp)$ is a field. In particular, $x\in\mp$. Hence $\Ass\,k(\mp)=\{\mp\}$. Furthermore, Lemma~\ref{commalg} implies that $\Supp\,k(\mp)=\overline{\Ass\,k(\mp)}=V(\mp)$. It is clear that $\Supp\,R/\mp=V(\mp)$, also see Exercise 3.19 in~\cite{Atiyah}.
\end{ex}






We recall the Nakayama's Lemma for convenience, which usually allows us to deduce contradictions when considering the annihilation properties of homologies.

\begin{lem}\label{nakayama}
(Nakayama) Let $R$ be a commutative ring with identity. Suppose $M$
is a finitely generated $R$-module and $\ma$ is an ideal of $R$
contained in the Jacobson radical of $R$. Then $\ma M=M$ implies
$M=0$.
\end{lem}

\pf See Proposition 2.6 in~\cite{Atiyah}.~~~\eop\\


The model category structure we are using here on the category $\Ch(R)$ of chain complexes comes from~\cite{Hovey}, in which the \textit{cofibrant} objects are exactly the \textit{DG-projective} complexes, i.e. \textit{K-projective} complexes such that each degree is projective. Here a complex $P$ is K-projective if and only if the functor $\Hhom_R(P,-)$ sends acyclic complexes to acyclic ones, or equivalently, $\Hom_{\D(R)}(P,X)\cong\Hom_{\K(R)}(P,X)$ for all $X$. See Proposition 1.4 in~\cite{Spalt}. One can consult the definition of \textit{total Hom chain complexes} $\Hhom_R(-,-)$ in Chapter V of~\cite{hilt} and properties of Hom functors in $\D(R)$ in Chapter 10 of~\cite{weib}

We denote by $k(\mp)$ the residue field of a prime ideal $\mp$, defined as $k(\mp)=R_\mp/\mp R_\mp$.

\begin{prop}\label{weakcpt}
Suppose $X\in\D^c_\fg(R)$. Then

\noindent(1) $\Hom(X,\Sigma^n\bigoplus_\a Y_\a)\cong\bigoplus\Hom(X,\Sigma^nY_\a)$ for all $R$-modules $Y_\a$;

\noindent(2) $X\otimes k(\mp)\cong\bigoplus_i\Sigma^i\bigoplus_{\a\in I(i)}k(\mp)$ for finite $I(i)$;

\noindent(3) $\Hom(X,Y)\otimes R_\mp\cong\Hom(X_\mp,Y_\mp)$ for every $Y\in\D^b(R)$;

\noindent(4) $X$ has a projective representative which is finitely generated in each degree.

\end{prop}

\pf For Property (1), since $X$ can be represented by a cofibrant object which is in particular K-projective, Proposition 1.4 in~\cite{Spalt} implies that it is sufficient to show the functor $\Hom(X_i,-)$ commutes with direct sums in the category of $R$-modules, which is clear since each $X_i$ is finitely generated. Property (2) comes from Lemma 2.17 in~\cite{Bokstedt}. Property (3) is part of our Proposition~\ref{passlocal}, and Property (4) is exactly the first part of Lemma 2.3.6 in~\cite{Hovey}, since the lifting property of cofibration implies the projectivity.~~~\eop\\





We end this section by showing that in $\D(R)$, under some sufficient conditions on the two variables of the Hom functor, with which the localization functor commutes. This is usually a key step when considering generalization from local to global. Nevertheless, this proposition is not always true in general.

\begin{prop}\label{passlocal}
Suppose one of the following conditions holds,

\noindent(1) $M\in \D_{\perf}(R)$ and $N\in \D(R)$,

\noindent(2) $M\in \D^+_\fg(R)$ and $N\in \D^-(R)$ or

\noindent(3) $M\in \D_\fg^c(R)$ and $N\in \D^b(R)$.

\noindent Then for every $\mp\in\Spec\,R$ we have an
isomorphism
\begin{displaymath}
\Hom_{\D(R)}(M,N)\otimes R_\mp\cong \Hom_{\D(R_\mp)}(M_\mp, N_\mp)
\end{displaymath}
as $R_\mp$-modules.

\end{prop}

\pf The key step is to observe that each degree $D_n=\prod_{p+q=n}\Hom_R(M_{-p},N_q)$ of the total Hom chain complex $\Hhom_R(M,N)$ becomes a coproduct. Take condition (1) for instance, other cases are similar except (3) which needs Property (4) of Proposition~\ref{weakcpt}. Since
$M$ is bounded, the product is finite, hence $D_n\cong
\bigoplus_{p+q=n}\Hom_R(M_{-p},N_q)$. Note that each $M_i$ is
finitely generated and $R_\mp$ is flat, $\Hom_R(M_{-p},N_q)\otimes
R_\mp\cong \Hom_{R_\mp}(M_{-p}\otimes R_\mp, N_q\otimes R_\mp)$ by
Proposition 2.10 in~\cite{eisen}. Thus $D_n\otimes R_\mp\cong
\bigoplus_{p+q=n} \Hom_{R_\mp}(M_{-p}\otimes R_\mp, N_q\otimes
R_\mp)$, which is the degree $n$ part of the chain complex
$\Hhom_R(M,N)\otimes R_\mp$, or $\Hhom_{R_\mp}(M_\mp,N_\mp)$ by the
definition of total Hom chain complex. Since $M$ is bounded below and projective in each degree, then we have
$\Hom_{\D(R)}(M,N)\cong H_0(\Hhom_R(M,N))$. Therefore,
\begin{eqnarray*}
\Hom_{D(R)}(M,N)\otimes R_\mp & \cong & H_0(\Hhom_R(M,N))\otimes
R_\mp
\cong H_0(\Hhom_R(M,N)\otimes R_\mp)\\
&\cong & H_0(\Hhom_{R_\mp}(M_\mp,N_\mp))\cong
\Hom_{D(R_\mp)}(M_\mp,N_\mp),
\end{eqnarray*}
where the second isomorphism comes from Proposition 2.5
in~\cite{eisen}.~~~\eop\\


\section{Aisles, $t$-structures and nullity classes}\label{atn}

In this section, we introduce the basic notions properties of nullity classes and related topics such as aisles and $t$-structures, and establish their fundamental relations. Denote by $\tl$ a triangulated category whose suspension functor is $\Sigma$. The basic properties of triangulated categories can be found in~\cite{BBD} and~\cite{weib}, for instance.

\begin{de}
A non-empty full subcategory $\al$ of $\tl$ is a preaisle if:

\noindent(1) for every $X\in \al$, we have $\Sigma X\in \al$;

\noindent(2) for every distinguished triangle $X\rightarrow
Y\rightarrow Z\rightarrow \Sigma X$, if $X,Z\in
  \al$ then $Y\in \al$.
\end{de}

A preaisle $\al$ is called \textit{cocomplete} if $\al$ is closed
under coproducts. A preaisle $\al$ is called an \textit{aisle} if
the inclusion functor $\al\hookrightarrow \tl$ admits a right
adjoint. As an easy consequence of the definition of preaisles, we obtain the following proposition.

\begin{prop}
Suppose $\al$ is a preaisle of $\tl$ then:

\noindent(1) $\al$ contains the zero object, i.e. $0\in \al$.

\noindent(2) if $X\in \al$ and $Y\in \tl$ is isomorphic to $X$, then
$Y\in \al$.

\end{prop}

\pf Since $\al\neq\emptyset$, we can take $X\in \al$. Then it
suffices to consider the following two distinguished triangles
$X\rightarrow 0\rightarrow \Sigma X\rightarrow \Sigma X$ and
$X\rightarrow Y\rightarrow 0\rightarrow\Sigma X$.~~~\eop\\

The concept of $t$-structures traces back to the work of A.A. Beilinson, J. Bernstein and P. Deligne in~\cite{BBD}, and that of aisles was introduced in the work of B. Keller and D. Vossieck in~\cite{Keller}. The Theorem~\ref{t-structure} implies that these two concepts are essentially the same thing.

\begin{de}
A $t$-structure on a triangulated category $\tl$ is a pair of two
subcategories $\al,\al'\subset \tl$ such that:

\noindent(1) $\Sigma \al\subset \al$ and $\al'\subset \Sigma\al'$;

\noindent(2) for every $A\in \al$ and every $B\in \Sigma^{-1}\al'$, we have
$\Hom(A,B)=0$;

\noindent(3) for every $X\in \tl$ there is a distinguished triangle $A\rightarrow
X\rightarrow B\rightarrow \Sigma A$ such that $A\in \al$ and $B\in
\Sigma^{-1}\al'$.
\end{de}

Let $\al\subset \tl$ be a subcategory. We define $\al^\perp=\{X\in
\tl~|~\Hom(A,X)=0$ for every $A\in \al\}$. Notice that $\al^\perp$
still lies in $\tl$. The following theorem is essentially the first proposition in~\cite{Keller}.

\begin{thm}\label{t-structure}
A preaisle $\al$ is an aisle, that is the inclusion
$\al\hookrightarrow \tl$ admits a right adjoint, if and only if
$(\al,\Sigma\al^\perp)$ is a $t$-structure on $\tl$.
\end{thm}

Associated to the $t$-structure $(\al,\Sigma\al^\perp)$, we denote
the right adjoint of the inclusion $\al\hookrightarrow \tl$ by
$(-)\langle \al\rangle:\tl\rightarrow \al$, and the left adjoint of
the inclusion $\al^\perp\hookrightarrow \tl$ by
$P_\al:\tl\rightarrow \al^\perp$. These adjoints are called
\textit{truncation functors}, and we often refer to $P_\al$ the
\textit{nullification functor} as well. Due to Theorem~\ref{t-structure}, for every
$X\in \tl$ we obtain a natural triangle,
\begin{displaymath}
\CD X\langle \al\rangle@>>> X@>>> P_\al(X)@>>> \Sigma X\langle
\al\rangle.
\endCD
\end{displaymath}

\noindent\textbf{Proof of Theorem~\ref{t-structure}.} Suppose
$(\al,\Sigma\al^\perp)$ is a $t$-structure and $X\in \tl$. Then we
have a distinguished triangle $A\rightarrow X\rightarrow
B\rightarrow\Sigma A$ with $A\in \al$ and $B\in\al^\perp$. Applying
the functor $\Hom_\tl(A',-)$ to this triangle for any $A'\in\al$ and
noticing that $\al\subseteq\tl$ is full, we have an exact sequence
\begin{displaymath}
\CD
\Hom_\tl(A',\Sigma^{-1}B)@>>>\Hom_\al(A',A)@>>>\Hom_\tl(A',X)@>>>\Hom_\tl(A',B),
\endCD
\end{displaymath}
where both ends are trivial since $B\in\al^\perp$ and $\al^\perp$ is
closed under desuspension by definition. Hence
$\Hom_\al(A',A)\rightarrow\Hom_\tl(A',X)$ is an isomorphism, which
is induced by $u:A=X\langle \al\rangle\rightarrow X$.

Conversely, suppose $\al$ is an aisle. The adjointness
$\Sigma\dashv\Sigma^{-1}$ of suspension functors implies $\al'=\Sigma
\al^\perp$ is closed under desuspension. Now suppose
\begin{displaymath}
\CD X\langle\al\rangle@>u>>X@>v>>Y@>w>>\Sigma X\langle\al\rangle
\endCD
\end{displaymath}
is a distinguished triangle, it remains to show that $Y\in
\al^\perp$. Indeed, let $A\in\al$ and $f\in\Hom_\tl(A,Y)$. Then we
have a morphism between two distinguished triangles by definition
\begin{displaymath}
\CD
X\langle\al\rangle@>h>>B@>>>A@>wf>>\Sigma X\langle\al\rangle\\
@|@VgVV @VfVV@|\\
X\langle\al\rangle@>u>>X@>v>>Y@>w>>\Sigma X\langle\al\rangle.
\endCD
\end{displaymath}
It follows from the adjunction $u_\ast$ and the commutativity of the
first square that $u=gh=u_\ast(g')h=ug'h$, i.e.
$g'h=id_{X\langle\al\rangle}$, for some $g':B\rightarrow
X\langle\al\rangle$. Thus $wf=0$. Hence $f$ factors through $X$ as
$f=vf'$ for some $f':A\rightarrow X$.
Therefore $f=vf'=vu_\ast(f'')=vuf''=0$ for some $f'':A\rightarrow
X\langle\al\rangle$, by using the adjunction $u_\ast$ again.~~~\eop\\

As an immediate consequence we have

\begin{coro}\label{nullify}
Suppose $\al$ is an aisle, then $P_\al(X)=0$ implies $X\in\al$.
\end{coro}

Next we move into the study of the derived category $\D(R)$ of a ring $R$ and its subcategories, nullity classes, basically following the work of Stanley, see his paper~\cite{Stan}.

\begin{de}
Let $\dl\subseteq \D(R)$ be a full triangulated subcategory and $\al$ be a cocomplete pre-aisle in $\D(R)$. A
nullity class in $\dl$ is a full subcategory of the form $\al\cap
\dl$.
\end{de}

we Let $E\in \D(R)$, and we denote by $\al=\overline{C(E)}$ the
smallest nullity class containing $E$ in $\D(R)$. We denote by
$P_E=P_\al$ the nullification functor. Let $F\in\D(R)$ be another object, we denote by $E<F$ if $P_E(F)=0$.


\begin{prop}\label{nullityclass}
Let $\dl\subset \D(R)$ be a full triangulated subcategory. Then every aisle in $\dl$ is a nullity class and every nullity class is a pre-aisle.
\end{prop}

Proposition~\ref{nullityclass} was proved as Proposition 2.14 in~\cite{Stan}, which tells us that nullity classes are slightly generalized aisles from preaisles.

\begin{prop}~\label{ntran}
Let $E,F\in \D(R)$. Then $E<F$ if and only if
$\overline{C(F)}\subseteq\overline{C(E)}$. Moreover,
$\overline{C(E)}=\{X\in\D(R)~|~P_E(X)=0\}$.
\end{prop}

\pf See Proposition 3.3 in~\cite{Stan}.~~~\eop\\

Clearly, Proposition~\ref{ntran} implies that the killing relation ``$<$'' is transitive. Also, it is easy to show by the Eilenberg Swindle that nullity classes are closed under retracts, since they are cocomplete and closed under suspensions and extensions.

\begin{lem}\label{retracts}
Let $\overline{C(E)}$ be a nullity class for some $E\in D(R)$. Then $\overline{C(E)}$ is closed under retracts.
\end{lem}

\pf See Lemma 3.2 in~\cite{Stan}.~~~\eop\\

Finally, we can deduce the following lemma from the advantage of being a nullity class.

\begin{lem}\label{nullityprop}
Suppose $E\in \D(R)$. If $\mp\notin \Supp H_i(E)$ for $i\leqslant
n$, then for every $X\in \overline{C(E)}$, we have $\mp\notin \Supp H_i(X)$
for $i\leqslant n$.
\end{lem}

\pf Notice that the condition $H_i(X)\otimes R_\mp=0$ for
$i\leqslant k$ is closed under suspensions, direct sums and
extensions. Indeed, the conclusion follows from the properties that
$H_i(\Sigma X)\cong H_{i-1}(X)$, the homology functors commute with
direct sums and the localization functor at $\mp$ is exact, see
Proposition 3.3 in~\cite{Atiyah}.~~~\eop

\section{Koszul complexes and annihilation}\label{kca}

As typical objects in the derived category $\D_\perf(R)$, the Koszul complexes play a crucial role in this paper. Their various properties can be found in~\cite{eisen} and~\cite{weib}, for example. We prove in this section as Lemma~\ref{genkilling} that in the light of the annihilation property of Koszul complexes, the quotients $R/\mp$ by $\mp$ which comes from the support of an object $M\in \D(R)$ are sufficient to generate $M$.

\begin{nota}\label{koszulcomp}
Let $\mp=(x_1,...,x_k)$ be an ideal in
$R$. Denote by $K(\mp)=K(x_1,...,x_k)$ the Koszul complex
associated to the sequence $x_1,...,x_k$ and
$K'=K(x_1^{n_1},...,x_k^{n_k})$ for some powers $n_1,...,n_k\in\NN$. Also we
define $K(0)=R$, and inductively define
\begin{eqnarray*}
K(i)=\Cone(K(i-1)\stackrel{x_i}{\rightarrow}
K(i-1)),
\end{eqnarray*}
as the mapping cone of the multiplication map by $x_i$. Moreover, for every $M\in \D(R)$ we define $M(i)=M\otimes K(i)$. In particular, $M(k)=M\otimes K(\mp)$. Similarly, we define $K'(0)=R$ and
$K'(i)=\Cone(K'(i-1)\stackrel{x_i^{n_i}}{\rightarrow} K'(i-1))$ for $n_i\in\NN$. In particular, $H_0(K(\mp))=R/\mp$, see Section 17.2 in~\cite{eisen}.
\end{nota}

Recall the annihilation property of Koszul complexes proved as Proposition 17.14 in~\cite{eisen}.

\begin{prop}\label{annihilation}
Suppose $M$ is an $R$-module and $\mp$ is an ideal of $R$ having $k$ generators. Then for
each $i=0,1,...,k$, the homology $H_i(M\otimes_R K(\mp))$ is annihilated by
every $x\in \mp$. In particular, such $x$ annihilates all the homologies
of $K(\mp)$.
\end{prop}

\begin{ex}\label{suppkosz}
Let $\mq=\ann^R(a)\in\Ass\,H_i(K(\mp))$ with $a\in H_i(K(\mp))$. Then Proposition~\ref{annihilation} implies that $a\cdot\mp=0$ thus $\mq\in V(\mp)$. Hence $\Supp\, K(\mp)\subseteq V(\mp)$. Also notice that $\Supp\,H_0(K(\mp))=V(\mp)$ by Example~\ref{suppkp}. Therefore, $\Supp\,K(\mp)=V(\mp)$.
\end{ex}

\begin{coro}\label{koszulhom}
Let $\mm$ be a maximal ideal of $R$ having $k$ generators. Then $H_i(K(\mm))\cong\bigoplus_{\Lambda_i}
k(\mm)$ for some index set $\Lambda_i$. In particular, $H_0(K(\mm))\cong H_k(K(\mm))\cong k(\mm)$.
\end{coro}

\pf Notice that the differentials of Koszul complex $K(\mm)$ is defined by multiplication by the generators of $\mm$, the computation follows immediately from Proposition~\ref{annihilation} since $\mm$ is maximal and $R/\mm\cong k(\mm)$.~~~\eop\\

The next technical result of annihilation roughly says that for a fixed element $x\in R$, every map with domain a perfect complex $X$ and arbitrary codomain $Y$ can be annihilated by some power of $x$ as soon as the homologies of $Y$ can be annihilated by some power of $x$. Therefore, in particular, it allows us to obtain a nonzero map from a Koszul complex to an appropriate object in $\D(R)$, as Lemma~\ref{minsupp} shows, which is crucial in the proof of Lemma~\ref{genkilling}.

\begin{lem}\label{homologyann}
Let $X\in\D_{\perf}(R)$ and $x\in R$. Let $Y\in \D(R)$ and $f:X\rightarrow Y$ be a morphism in $\D(R)$. Suppose for every $\ast$ and every $\a\in H_\ast(Y)$, there exists $n\in\NN$ such that $x^n\a=0$. Then $x^lf=0$ for some $l\in \NN$.
\end{lem}



\pf This is proved by induction using the filtration
$X(i)=\bigoplus_{j=0}^iX_j$ on $X$, assuming that $X$ is bounded
between dimension 0 and $k\geqslant0$ with each degree free on
finitely many generators.

First we show the case for $k=0$, i.e. $X$ has nontrivial homology
for only degree zero. Then that $X$ is perfect implies
$\Hom_{\D(R)}(X,Y)\cong\bigoplus\Hom_{\D(R)}(R,Y)\cong \bigoplus
H_0(Y)$, which has only finitely many summands by assumption. Hence
any $f\in \Hom_{\D(R)}(X,Y)$ is annihilated by some power of $x$.



In general, for $k\geq 1$ we set $f^{(k-1)}_{k-1}=f^{(k-1)}i_{k-1},
f^{(k-1)}_k=f^{(k-1)}i_k$ and $f^{(k-1)}=x^{n_{k-2}}f^{(k-2)}$ such
that the following diagram commutes
\begin{displaymath}
\CD X(k-1) @>i_{k-1,k}>> X(k) @>i_k>> X\\
@Vf^{(k-1)}_{k-1}VV  @Vf^{(k-1)}_kVV   @VVf^{(k-1)}V\\
Y @= Y @= Y
\endCD
\end{displaymath}
where $i_{k-1}=i_k\circ i_{k-1,k}$ and $f^{(0)}=f$. We replace the three vertical
maps by $f^{(k)}_{k-1}=x^{n_{k-1}}f^{(k-1)}_{k-1}=0,
f^{(k)}_k=x^{n_{k-1}}f^{(k-1)}_k$ and
$f^{(k)}=x^{n_{k-1}}f^{(k-1)}$, still preserving the commutativity, where $f^{(k)}_{k-1}=0$ for some power $n_{k-1}$ by induction.
Hence there exists a map $\rho_k$ making
the following diagram commute
\begin{displaymath}
\CD X(k-1) @>i_{k-1,k}>> X(k) @>\pi_k>> X(k)/X(k-1) @>>>\Sigma X(k-1)\\
@Vf^{(k)}_{k-1}=0VV  @Vf^{(k)}_kVV @V\exists V\rho_kV  \\
Y @= Y @= Y
\endCD
\end{displaymath}
so that $x^{n_k}\rho_k=0$ for some $n_k$ implies
$x^{n_k}f^{(k)}_k=0$. Therefore, since $X(k)=X$ thus $i_k=id_X$ and
$f^{(k-1)}=x^{\sum_{i=0}^{k-2}n_i}f$ is defined recursively, we
deduce that
\begin{eqnarray*}
0=x^{n_k}f^{(k)}_k=x^{n_k}x^{n_{k-1}}f^{(k-1)}_k=x^{n_k}x^{n_{k-1}}f^{(k-1)}i_k=x^{\sum^k_{i=0}n_i}f,
\end{eqnarray*}
as we required.~~~\eop

\begin{lem}\label{minsupp}
Let $M\in D(R)$. Suppose $\mp=(x_1,...,x_k)\in\Supp M$ is minimal such that $\mp\in \Supp H_n(M)$ for some $n$. Then $\Hom_{D(R)}(\Sigma^nK'(\mp),M)\neq0$, for appropriate powers.
\end{lem}

\pf Notice that $K'(\mp)$ is a perfect complex, by Proposition~\ref{passlocal} it suffices to show that $\Hom_{D(R_\mp)}(\Sigma^nK'(\mp)_\mp,M_\mp)\neq0$. Indeed, for every $\ast$ for which $H_\ast(M)\neq 0$ and, for a decomposition of a cyclic module generated by a nontrivial element $\a\in H_\ast(M)$, we can deduce that
\begin{displaymath}
\mp_i R_\mp\in\Ass R_\mp/\mp_i R_\mp\subseteq\Ass (R\a)_\mp\subseteq\Ass H_\ast(M_\mp)\subseteq V(\mp R_\mp),
\end{displaymath}
where the last inclusion holds because the localization at $\mp$ kills all the primes containing $\mp$ and the minimality of $\mp\in\Supp M$ implies that all such primes $\mp_iR_\mp$ must be identified with $\mp R_\mp$. Inductively for each $j=1,...,k$ there is some $n_j$ such that $\a$ is annihilated by $x_j^{n_j}$.

For simplicity, we now assume that $R$ is local. Suppose $0\neq f_0\in \Hom(K'(0),M)\cong H_0(M)$. Since for any $\a\in
H_\ast(M)$ some power of $x_1$ annihilates $\a$, there is by
Lemma~\ref{annihilation} some $n_1\in\NN$ such that $x_1^{n_1}$
annihilates $f_0$, hence the map $f_0$
extends to a map $f_1:K'(1)\rightarrow M$ by considering the
distinguished triangle
\begin{displaymath}
\CD K'(0)@>x_1^{n_1}>> K'(0) @>>> K'(1) @>>> \Sigma K'(0).
\endCD
\end{displaymath}
Observe that the fact $H_0(f_0)\neq0$ implies that $H_0(f_1)\neq 0$
since $f_0$ factors through $f_1$. Now suppose the maps
$f_j:K'(j)\rightarrow M$ are constructed such that $H_0(f_j)\neq0$
for $j\leqslant k-1$. By Lemma~\ref{annihilation} again, that some
power of $x_k$ annihilates $H_\ast(M)$ implies the existence of some
$n_k\in\NN$ such that $x_k^{n_k}f_{k-1}=0$. Now consider the
distinguished triangle
\begin{displaymath}
\CD K'(k-1)@>x_k^{n_k}>> K'(k-1) @>>> K'(k)=K' @>>> \Sigma K'(k-1),
\endCD
\end{displaymath}
the annihilation implies there is an extension $f=f_k:K'\rightarrow M$. Our inductive assumption implies that
$H_0(f)\neq0$, as we required.~~~\eop\\

In particular, if $M\in\D(R)$ satisfies $\Supp M=V(\mp)=\Supp H_n(M)$ for some prime $\mp$, then Lemma~\ref{minsupp} applies. We end this section by proving a generalization of Lemma 5.3 in~\cite{Stan}, which shows the same thing but in the case when the object $M$ is bounded.

\begin{lem}\label{genkilling}
Let $M\in D(R)$. Then $\bigoplus_{i\in\ZZ}\bigoplus_{\mp\in \Supp H_i(M)}\Sigma^i R/\mp<M$.
\end{lem}

\pf Denote by $E=E(M)=\bigoplus_{i\in\ZZ}\bigoplus_{\mp\in \Supp H_i(M)}\Sigma^i R/\mp$. We prove it by contradiction, assuming that $P_EM\neq 0$. Notice that $E<M\langle E\rangle$ so that $\Ass H_{i-1}M\langle E\rangle\subseteq\Supp H_{i-1}(E)\subseteq \Supp H_{i-1}(M)$. Thus $\Ass H_i P_EM\subseteq \Supp H_i(M)\cup \Ass H_{i-1}M\langle E\rangle\subseteq\Supp H_i(M)\cup\Supp H_{i-1}(M)$ by Lemma 5.4 in~\cite{Stan}. Since every Noetherian local ring has finite dimension (see Corollary 11.11 in~\cite{Atiyah}), one can pick up a minimal prime in $\bigcup_i \Ass H_iP_EM$ such that $\mp\in\Supp H_n(M)$ for some $n$. The natural distinguished triangle
\begin{displaymath}
M\langle E\rangle\rightarrow M\rightarrow P_EM\rightarrow \Sigma M\langle E\rangle
\end{displaymath}
allows us to identity $\mp$ to be minimal also in $\Supp M$. Hence by Lemma~\ref{minsupp} there is a map $f:\Sigma^nK(\mp)\rightarrow P_EM$ such that $H_n(f)\neq 0$. Therefore, the fact $E<\Sigma^nK(\mp)$ by Lemma 5.3 in~\cite{Stan} implies $\Hom(\Sigma^jE,P_EM)\neq0$ for some $i\geq0$ by Proposition 3.4 in~\cite{Stan}, a contradiction.~~~\eop

\section{The key lemma}\label{key}

This section presents the key ingredient, Lemma~\ref{crucial} , in order to produce an invariant of nullity classes. In fact, this lemma implies the injectivity of the map $\phi$ defined in the next section, which maps the class of aisles in the derived category $\D^c_\fg(R)$ into the class of perversity functions on the power set $\pl(\Spec R)$.

We first recall two facts which are Lemma 6.1 and Lemma 6.2 respectively in~\cite{Stan}.

\begin{lem}\label{k1}
If $E<F$, then $E\otimes M< F\otimes M$ for every $M\in \D(R)$.
\end{lem}

\begin{lem}\label{k2}
If $M\in\D(R)$ such that $H_i(M)=0$ for $i<0$, then $R<M$.
\end{lem}


In Notation~\ref{koszulcomp} we define the Koszul complex $K(i)$ as the mapping cone of $x_i:K(i-1)\rightarrow K(i-1)$ with $K(0)=R$ and $M(i)=K(i)\otimes M$.

\begin{lem}\label{localcase}
Let $M\in\D_\fg(R)$ and $H_n(M)\neq 0$ for some $n$. Then

\noindent(1) if $(R,m)$ is local, then $H_n(M\otimes K(\mm))\neq0$ and,

\noindent(2) in general, there is a $\mp\in\Supp H_n(M)$ such that $\mp\in\Supp H_n(M\otimes K(\mp))$.
\end{lem}

\pf Suppose $R$ is local with $\mm=(x_1,...,x_k)$. Then we have distinguished triangles
\begin{displaymath}
\CD M(i)@>x_{i+1}>>M(i)@>>>M(i+1)
\endCD
\end{displaymath}
which induce long exact sequences on homologies
\begin{displaymath}
\CD H_n(M(i))@>x_{i+1}>>H_n(M(i))@>>>H_n(M(i+1)).
\endCD
\end{displaymath}
It turns out inductively on $i$, the homology
$H_n(M(i+1))\neq0$ because otherwise, the induced map $x_{i+1}$ on
the finitely generated module $H_n(M(i))$ is surjective, so that Lemma~\ref{nakayama} implies $H_n(M(i))=0$, contradicting our assumption.

In general, since $H_n(M)\neq0$, there is a prime ideal $\mathfrak{p}$ of
$R$ such that $H_n(M_\mp)\cong H_n(M)\otimes
R_\mp\neq 0$ by Proposition 3.8 in~\cite{Atiyah}. Then
$H_n(M\otimes K(\mp))\otimes R_\mp\cong
H_n(M_\mp\otimes K(\mp)_\mp)\cong
H_n(M_\mp\otimes K(\mp R_\mp))\neq0$ by (1) since
$(R_\mp,\mp R_\mp)$ is local and
$M_\mp\in D_\fg(R_\mp)$.~~~\eop\\


It is well-known that associated to each standard truncation $\t_{\geq n}$ of an object $M$ in $\D(R)$ there is a distinguished triangle
\begin{eqnarray*}
\t_{\geq n}M\rightarrow M\rightarrow\t_{< n}M
\end{eqnarray*}
and also that $\t_{\geq n}\t_{\leq n}M\cong\t_{\leq n}\t_{\geq n}M\cong \Sigma^nH_n(M)$, where $\t_{\leq n}=\t_{<n+1}$. In particular, we have a distinguished triangle
\begin{eqnarray*}
\t_{\geq n+1}M\rightarrow \t_{\geq n}M\rightarrow\Sigma^nH_n(M)
\end{eqnarray*}
for each $n$. See for example~\cite{BBD} and~\cite{weib}.

\begin{lem}\label{transit2}
Let $(R,\mm)$ be a local ring whose maximal ideal has $k$ generators and $M\in \D_\fg(R)$. Suppose $H_i(M\otimes k(\mm))=0$ for
$n-k\leq i\leq n$. Then $H_n(M\otimes K(\mm))=0$.



\end{lem}

\pf Truncating $K(\mm)$ each step a homology $\Sigma^iH_i(K(\mm))$, we have a filtration
\begin{displaymath}
0=K_{k+1}\rightarrow K_k\rightarrow K_{k-1}\rightarrow\cdots\rightarrow K_1\rightarrow K_0=K(\mm)
\end{displaymath}
for $K(\mm)$ with $K_i=\tau_{\geq i}K(\mm)$ such that $K_{i+1}\rightarrow K_i\rightarrow \Sigma^iH_i(K(\mm))$ is a distinguished triangle for each $i=0,1,...,k$. Since tensor product preserves distinguished triangles, the sequence
\begin{displaymath}
0=M\otimes K_{k+1}\rightarrow M\otimes K_k\rightarrow M\otimes K_{k-1}\rightarrow\cdots\rightarrow M\otimes K_1\rightarrow M\otimes K_0=M\otimes K(\mm)
\end{displaymath}
remains a collection of distinguished triangles. Hence by taking the $n$-th homology and the long exact sequences, $H_i(M\otimes k(\mm))=0$ for $n-k\leq i\leq n$ implies $H_n(M\otimes \Sigma^iH_i(K(\mm)))\cong H_{n-i}(M\otimes (\oplus k(\mm)))=0$ for $i=0,1,...,k$ by Corollary~\ref{koszulhom}. Consequently, $H_n(M\otimes K_i)=0$ for each $i$. In particular, $H_n(M\otimes K(\mm))=0$ as required.~~~\eop


\begin{coro}\label{transit}
Suppose $M\in \D_\fg(R)$ and $\mp$ is a prime ideal having $k$ generators, such that $H_i(M\otimes k(\mp))=0$ for $n-k\leq i\leq n$. Then $H_n(M_\mp)=0$.




\end{coro}

\pf Notice that $M\otimes k(\mp)\cong M_\mp\otimes k(\mp R_\mp)$. Thus $H_n(M_\mp\otimes K(\mp R_\mp))=0$ by Lemma~\ref{transit2}. Hence $H_n(M_\mp)=0$ thanks to Lemma~\ref{localcase}.~~~\eop\\


The following lemma is a generalization of Lemma 6.4 in~\cite{Stan} in which the object is assumed to be bounded.

\begin{lem}\label{killkp}
Suppose $M\in D_\fg(R)$ and $\mp\in\Supp H_n(M)$ has $k$ generators. Then $M<\Sigma^{n-i} k(\mp)$ for some $0\leq i\leq k$. In particular, $M<\Sigma^{n-i} k(\mq)$ for every $\mq\in V(\mp)$.
 \end{lem}

\pf Observe that $R<k(\mp)$ implies $M<M\otimes k(\mp)$ by Lemma~\ref{k1}
and~\ref{k2}, which is a direct sum of suspensions of finite copies of $k(\mp)$ by
Lemma 2.17 in~\cite{Bokstedt}. Thus $\mp\in \Supp H_n(M)$ implies $H_{n-i}(M\otimes k(\mp))\neq0$ for some $0\leq i\leq k$ by Lemma~\ref{transit}. Hence $M<\Sigma^{n-i}k(\mp)$, since nullity class is closed under retracts. The second statement follows from the fact that the support sets are specialization closed.~~~\eop\\



The proof of the next lemma follows precisely that of Proposition 6.6 in~\cite{Stan}. For convenience, we denote by $\bar{k}(\mp)=\bigoplus_{\mq\in V(\mp)}k(\mq)$.

\begin{lem}\label{kil}
Suppose $\dim R/\mp$ is finite. Then $\bar{k}(\mp)<\Sigma^{\dim R/\mq} R/\mq$ for every $\mq\in V(\mp)$.
\end{lem}

\pf Prove by induction on $\dim R/\mp$. Suppose $\dim R/\mp=0$, i.e. $\mp=\mm$ is maximal, then the inequality holds since $k(\mm)\cong R/\mm$. Now for a fixed $\mq\in V(\mp)$, assume $\bar{k}(\mp)<\Sigma^{\dim R/\mq'}R/\mq'$ for every $\mq'\in V(\mq)-\{\mq\}$. Consider the short exact sequence of $R$-modules
\begin{displaymath}
0\rightarrow R/\mq\stackrel{f}{\rightarrow} k(\mq)\rightarrow M\rightarrow 0
\end{displaymath}
where $f$ is the natural map and $M$ its quotient. Since $f\otimes R_\mq$ becomes an isomorphism, $\Supp M\subseteq V(\mq)-\{\mq\}$. Hence by induction together with Lemma~\ref{genkilling}, we can deduce that $\bar{k}(\mp)<\Sigma^{\dim R/\mq-1}\bigoplus_{\mq'\in\Supp M} R/\mq'<\Sigma^{\dim R/\mq-1}M$. Furthermore, the distinguished triangle
\begin{displaymath}
\Sigma^{\dim R/\mq-1}M\rightarrow \Sigma^{\dim R/\mq}R/\mq\rightarrow\Sigma^{\dim R/\mq}k(\mq)
\end{displaymath}
allows us to conclude that $\bar{k}(\mp)<\Sigma^{\dim R/\mq}R/\mq$, since nullity classes are closed under retracts and extensions.~~~\eop\\


As an easy consequence we can show the following statement.

\begin{coro}\label{ringdim}
Suppose $\dim R/\mp$ is finite. Then $\bar{k}(\mp)<\Sigma^{\dim R/\mp} K'(\mp)$ for every $K'$ with arbitrary finite powers on the generators of $\mp$. In particular, $k(\mm)<K'(\mm)$ if $\mm$ is maximal.
\end{coro}

\pf By Lemma~\ref{genkilling}, $\bigoplus_{i\geq0}\bigoplus_{\mq\in \Supp H_i(K'(\mp))}\Sigma^iR/\mq<K'(\mp)$, thus by Lemma~\ref{kil} we have $\bar{k}(\mp)<\Sigma^{\dim R/\mp}\bigoplus_{i\geq 0}\bigoplus_{\mq\in \Supp H_i(K'(\mp))} \Sigma^iR/\mq<\Sigma^{\dim R/\mp}K'(\mp)$, since nullity class is closed under coproducts and suspensions.~~~\eop\\

As a crucial step in proving our key lemma, we show that

\begin{lem}\label{thereismap}
Suppose $M\in D_\fg^c(R)$ and $\mp\in\Supp H_n(M)$, $\mp\notin\Supp H_i(M)$ for $i<n$. Then there is a map $f:M\rightarrow \Sigma^nH_n(M)$ such that $H_n(f_\mp)$ is an isomorphism.
\end{lem}

\pf Notice that $H_i((\t_{<n}M)_\mp)=H_i(\t_{<n}M)_\mp=0$ for all $i<n$. Thus the natural map $(\t_{\geq n}M)_\mp\rightarrow M_\mp$ associated to the standard truncation of $M$ at $n$ then localized at $\mp$ is a quasi-isomorphism. Hence
\begin{displaymath}
\Hom((\t_{\geq n}M)_\mp,\Sigma^nH_n(M)_\mp)=\Hom(M_\mp,\Sigma^nH_n(M)_\mp)
=\Hom(M,\Sigma^nH_n(M))\otimes R_\mp
\end{displaymath}
by Proposition~\ref{passlocal}. Therefore, there is a morphism $f:M\rightarrow \Sigma^n H_n(M)$ corresponding to the natural map $\pi:\t_{\geq n}M\rightarrow\Sigma^nH_n(M)$ via the above isomorphisms. Since $H_n(\pi_\mp)$ is an isomorphism, so is $H_n(f_\mp)$.~~~\eop\\





Finally we are well-prepared to prove the key lemma.

\begin{lem}\label{crucial}
Suppose $\dim R$ is finite and $M\in \D^c_\fg(R)$. Then $M<\Sigma^n R/\mp$ for every $\mp\in \Supp H_n(M)$.
\end{lem}


\pf It is not hard to show that for every degree $m$ where $\mp$ stands, we have $M<\Sigma^{m+\dim R/\mp}R/\mp$ by Lemma~\ref{killkp} and Lemma~\ref{kil} together with the transitivity of nullity classes. Define
\begin{displaymath}
l(\mp)=\inf\{i\in\ZZ~|~\mp\in \Supp H_i(M)\}.
\end{displaymath}
If $l(\mp)=-\infty$, then for each degree $n$ where $\mp$ stands, $M<\Sigma^{m+\dim R/\mp}R/\mp<\Sigma^nR/\mp$ since $m$ can be chosen such that $n-m-\dim R/\mp\geq0$. So fix a prime $\mp$ and suppose $l=l(\mp)$ is finite, i.e. $\mp\in\Supp H_l(M)$ and $\mp\notin\Supp H_i(M)$ for $i<l$.  Also, suppose for every $n$ and every $\mq\in\Supp H_n(M)$, $M<\Sigma^{n+k}R/\mq$ for $1\leq k\leq \dim R/\mq$. We need to show that $M<\Sigma^lR/\mp$ by downward induction.

The finiteness of $l(\mp)$ implies that there is a map $f:M\rightarrow \Sigma^lH_l(M)$ by Lemma~\ref{thereismap}. Now consider the distinguished triangle
\begin{displaymath}
M\stackrel{f}{\rightarrow} \Sigma^lH_l(M)\rightarrow M'
\end{displaymath}
and the associated long exact sequence which splits into short ones as
\begin{displaymath}
0\rightarrow H_{l+1}(M')\rightarrow H_l(M)\stackrel{H_l(f)}{\longrightarrow} H_l(M)\rightarrow H_l(M')\rightarrow H_{l-1}(M)\rightarrow0
\end{displaymath}
and $H_{l+i+1}(M')\cong H_{l+i}(M)$ for either $i\geq1$ or $i\leq-2$. Since $H_l(f_\mp)$ is an isomorphism, thus $H_l(M')_\mp=0$ so that $\Supp H_l(M')\subseteq \Supp H_l(M)\cup\Supp H_{l-1}(M)-\{\mp\}$ by Lemma~\ref{commalg}. Next we show that $M<\Sigma^lR/\mq$ for every $\mq\in V(\mp)-\{\mp\}$. Since $l$ is the minimal number at which $\mp$ stands, Lemma 6.7 in~\cite{Stan} implies that $M<N=M\otimes K(\mq)$ such that $\mq\in\Supp H_l(N)$ and $\Supp N\subseteq V(\mq)$. Therefore, it is sufficient to show $N<\Sigma^lR/\mq$ for every $\mq\in V(\mp)-\{\mp\}$, which is proved by induction on the dimension $\dim R/\mq$, similarly to the proof of Lemma~\ref{kil}.

In fact, if $\mq=\mm$ is maximal then $N<\Sigma^l R/\mm=\Sigma^lk(\mm)$ by Lemma~\ref{killkp}. Now suppose it is true for every $\mq\in V(\mp)$ such that $\dim R/\mq\leq\dim R/\mp-2$, we need to show it is true for every $\mq\in V(\mp)-\{\mp\}$, i.e. $\dim R/\mq\leq\dim R/\mp-1$. Consider the short exact sequence of $R$-modules
\begin{displaymath}
0\rightarrow R/\mq\stackrel{\iota}{\rightarrow} k(\mq)\rightarrow A\rightarrow 0
\end{displaymath}
where $\iota$ is the natural map and $A$ its quotient. Since $\iota_\mq$ becomes an isomorphism, $\Supp A\subseteq V(\mq)-\{\mq\}$. Hence by induction together with Lemma~\ref{genkilling}, we can deduce that $N<\Sigma^l\bigoplus_{\mq'\in\Supp A} R/\mq'<\Sigma^lA$. Furthermore, $N<\Sigma^lk(\mq)$ by Lemma~\ref{killkp}. Hence the distinguished triangle
\begin{displaymath}
\Sigma^lR/\mq\rightarrow\Sigma^lk(\mq)\rightarrow \Sigma^lA
\end{displaymath}
allows us to conclude that $N<\Sigma^lR/\mq$ since nullity classes are closed under extensions.

Combine what we have shown that $M<\Sigma^lR/\mq$ for every $\mq\in V(\mp)-\{\mp\}$ and the fact that the support of $M'$ is shifted down by one into the support of $M$ for every degree, we deduce by the downward induction and Lemma~\ref{genkilling} that $M<\bigoplus_i\bigoplus_{\mp\in\Supp H_i(M')}\Sigma^iR/\mp<M'$. Hence $M<\Sigma^lH_l(M)$ since nullity classes are closed under extensions. Therefore, $M<\Sigma^lR/\mp$ by Lemma 6.8 in~\cite{Stan}, and the proof is complete.~~~\eop\\

Therefore, together with Lemma~\ref{genkilling} our key lemma implies that in fact, each object in $\D(R)$ is completely determined by its support in the sense of construction of nullity classes, i.e. they generate the same category.

\section{A complete invariant}

Now we are ready to produce an invariant of nullity classes as well as aisles in the derived category $\D^c_\fg(R)$, as soon as we establish the involved maps $\phi$ and $N$ as follows.\\

We call a function $f:\ZZ\rightarrow \pl(\Spec\,R)$ a
\textit{perversity function} if $f(n)=\overline{f(n)}$, i.e. specialization closed, and $f(n)\subseteq f(n+1)$ for each $n\in \ZZ$. Then for every nullity
class $\al$ in a full triangulated subcategory $\dl\subseteq \D(R)$
we define $\phi(\al)$ by setting $\phi(\al)(n)=\{\mp\in
\Spec\,R~|~\mp\in \Supp H_n(M)$ for some $M\in\al\}$. Clearly such a
function is increasing since nullity class is closed under
suspensions and also $\phi(\al)(n)$ is closed under specialization
since the support set is. Hence the map
\begin{displaymath}
\phi:\{\text{nullity classes in $\dl$}\}\rightarrow
\{\text{perversity functions}\}
\end{displaymath}
is well-defined. On the other hand, let $S(f)=\bigoplus_{n\in\ZZ}\bigoplus_{\mp\in f(n)}\Sigma^nR/\mp$ and
$\overline{C(S(f))}$ be the smallest nullity class containing
$S(f)$. We define a function
\begin{displaymath}
N:\{\text{perversity functions}\}\rightarrow \{\text{nullity classes
in $\dl$}\}
\end{displaymath}
which associates each perversity function $f$ a nullity class
$N(f)=\overline{C(S(f))}\cap \dl$.

\begin{lem}
The maps $N$ and $\phi$ are order preserving.
\end{lem}

\pf Clearly $\phi$ is order preserving by definition. For $N$ we
consider $f,g$ any two perversity functions such that $f\leqslant
g$, i.e. $f(n)\subseteq g(n)$ for all $n$. Then $S(f)$ is a retract
of $S(g)$ by definition, which implies that
$S(f)\in\overline{C(S(g))}$ hence
$\overline{C(S(f))}\subseteq\overline{C(S(g))}$.~~~\eop

\begin{thm}~\label{mainthm}
Suppose $\dim R$ is finite and $\dl=\D^c_\fg(R)$. Then $\phi:\{$nullity classes in $\dl\}\rightarrow\{$perversity functions$\}$ is a bijection with inverse given by $N$.
\end{thm}

\pf Let $\al$ be a nullity class in $\dl$. Take $M\in\al$. Then by Lemma~\ref{retracts} we have $S(\phi(\al))=\bigoplus_{n\in\ZZ}\bigoplus_{\mp\in \phi(\al)(n)}\Sigma^nR/\mp<\bigoplus_{n\in\ZZ}\bigoplus_{\mp\in \Supp H_n(M)}\Sigma^nR/\mp$ since the latter coproduct is a retract of the former one. Thus $S(\phi(\al))<M$ by Lemma~\ref{genkilling}, so that $M\in\overline{C(S(\phi(\al)))}$. Hence $M\in N\phi(\al)$. Conversely, for every $\mp\in\phi(\al)(n)$ there is $M(\mp,n)\in\al$ such that $\mp\in\Supp H_n(M(\mp,n))$ by definition. Thus by Lemma~\ref{crucial} we have $M(\mp,n)<\Sigma^nR/\mp$, so that $M=\bigoplus_{n\in\ZZ}\bigoplus_{\mp\in \phi(\al)(n)}M(\mp,n)<\bigoplus_{n\in\ZZ}\bigoplus_{\mp\in \phi(\al)(n)}\Sigma^nR/\mp=\overline{S(\phi(\al))}$. Hence $\overline{S(\phi(\al))}\in\overline{C(M)}\subseteq\al$.

On the other hand, we show that $\phi N(f)(n)=f(n)$ for every perversity function $f:\ZZ\rightarrow \pl(\Spec R)$ and every $n\in\ZZ$. Take $\mp\in f(n)$. Then $S(f)<\Sigma^nR/\mp$ since $\Sigma^nR/\mp$ is a retract of $S(f)$, so that $\Sigma^nR/\mp\in N(f)$. It follows that $\mp\in\Ass H_n(\Sigma^nR/\mp)=\Ass R/\mp$, thus $\mp\in\phi(N(f))(n)$. Conversely, suppose $\mp\in\phi(N(f))(n)$. Then there is some $M\in N(f)=\overline{(S(f))}$ such that $\mp\in\Supp H_n(M)$. Lemma~\ref{nullityprop} implies that $\mp\in\Supp H_i(S(f))$ for some $i\leq n$. Hence $\mp\in f(i)\subseteq f(n)$, as required.~~~\eop\\

This theorem has immediate consequences as the following corollaries.

\begin{coro}
Every aisle in $\dl=\D_\fg^c(R)$ with finite $\dim R$ is of the form $\overline{C(E)}\cap
\dl$ for some $E\in \D(R)$.
\end{coro}

\begin{coro}
Let $\dl=\D_\fg^c(R)$ with finite $\dim R$. Suppose $\al$ and $\al'$ are aisles in $\dl$. Then
$\al\subseteq \al'$ if and only if $\phi(\al)\leqslant\phi(\al')$. Thus $\phi:\{$aisles in $\dl\}\rightarrow\{$perversity functions$\}$
is injective.
\end{coro}

Since aisles are nullity classes, the map $\phi$ gives a complete invariant for aisles (or $t$-structures) thanks to Theorem~\ref{mainthm}.

\bibliography{mybib}
\bibliographystyle{plain}

\end{document}